\documentclass
[12pt,a4paper]{article}
\usepackage[cp1251]{inputenc}
\usepackage{amssymb, amsfonts, amsmath, latexsym}
\usepackage[russian]{babel}
\begin{document}

\Large
\centerline{\bf Generating countable groups by discrete subsets}\vspace{6 mm}

\normalsize\centerline{\bf  Igor Protasov}\vspace{6 mm}

{\bf Abstract.} Every countable topological group $G$ has a closed discrete subset $A$ such that $G=AA^{-1}.$


\vspace{6 mm}

\centerline{\bf 1.	Introduction}

\vspace{6 mm}

This note is to strengthen radically the following statement [1, Theorem 2.2]: every countable topological group can be generated by some closed discrete subset. All topologies under   consideration  are supposed to be Hausdorff.
\vspace{3 mm}

{\bf Theorem.} {\it Every countable topological group $G$ has a closed discrete subset $A$ such that $G=AA^{-1}.$
\vspace{3 mm}

Proof.}We enumerate $G=\{ g_{n} : n\in \omega\}$, $g_{0}=e$,  $e$ is the identity of $G$,  and split the proof into two cases: $G$ is precompact (Case 1) and $G$ is not precompact  (Case 2). We use two equivalent definitions of precompact groups: $G$ is a subgroup of some compact group and, for every neighborhood $U$  of $e,$ there exists  a finite subset $F$  of $G$ such that $G=FU$.
\vspace{3 mm}

{\it Case 1.}
Let $G$ be a dense subgroup of a compact group $H$ and let $\mu$ be the Haar measure on $H$. We choose a sequence $(r_{n})_{n\in\omega} $ of positive real numbers such that $\Sigma_{i\in\omega} r_{i} < \frac{1}{6} .$

We put $x_{0}=y_{0}=e$, choose closed in $H$  neighborhood $U$  of $e,$ such that $\mu(U)< r_{0}$, denote $U_{0} = U$  and assume that, for some $n\in\omega $, we have chosen elements $x_{0} , …, x_{n}$  and $y_{0}, …, y_{n}$ of $G$, and closed in  $H$ neighborhoods $U_{0}, …, U_{n}$  of $e,$  such that , for each $k\leq  n$  and  $A_{k}= \{ x_{i}, y_{i}: i\leq  n\}$,

(1)	 $g_{k} \in A_{k} A_{k} ^{-1}$;

(2)	the subsets $\{x_{i} U_{i} : i\leq  k\}$  are pairwise disjoint, the subsets  $\{y_{i} U_{i}: i\leq k\}$ are pairwise disjoint and $ x_{i} U_{i}\bigcap y_{j}U_{j}=\emptyset$  for each $i,j\in \{1, ..., n\}$;

(3)	 $\mu(U_{i})< r_{i}$  for each $i\leq  k$.
\vspace{3 mm}

We suppose also that, for $n>1$, there is a numeration $z_{0}, …, z _{m(n)}$ of $A_{n}  A_{n} ^{-1}\setminus   A_{n}$ and closed in $H$ neighborhoods $V_{0}, …, V_{m(n)}$  of $e$  such that , for each $k\in  \{1, …, n\}$,     $A_{k}A_{k}^{-1}\setminus \{A_{k}\bigcup A_{k-1}A_{k-1}^{-1}\}=\{z_{m(k-1)}, ...,   z_{m(k)}\}$, and

(4)	$z_{i} V_{i}\bigcap A_{k}=\emptyset$,  $i\in\{0, ..., m(k)\};$

(5)	 $\mu(v_{i})< r_{i}, $  $i\in\{0, ..., m(k)\}.$

\vspace{3 mm}

To make the inductive step from $n$ to  $n+1$, we take the first element $g\in\{g_{i}: i\in\omega\}\setminus A_{n} A_{n} ^{-1}$ and denote

$$B=\{x_{i} U_{i} : i\leq n\}\bigcup \{y_{i} U_{i} : i\leq n\}\bigcup \{z_{i} V_{i} : i\leq m(n)\},  \  \   C=H\setminus B .$$

By (3), (5)  and the choice of $(r_{n})_{n\in\omega}$, we have $\mu(C)>\frac{1}{2}$. Hence,  $gC \bigcap C\neq\emptyset$.   Since $C$ is open  in $H$,  there are $x,y \in G \bigcap C$  such that $x=gy$.  We put $x _{n+1}= x$,  $y _{n+1}=y$, $A_{n+1}=A_{n}\bigcup \{x_{n+1}, y_{n+1}\}$ and note that $g\in A_{n+1}A_{n+1}^{-1}$. Since  $x _{n+1}, y _{n+1}\in G\setminus B$, there is a closed in $H$ neighborhood $U_{n+1}$ of $e$  such that $x_{n+1}U_{n+1}\bigcap B=\emptyset$, $y_{n+1}U_{n+1}\bigcap B=\emptyset$, $x_{n+1}U_{n+1}\bigcap y_{n+1} U_{n+1}=\emptyset$ and $\mu(U_{n+1})< r_{n+1}$. Hence,  (1), (2) (3)  are satisfied for $k=n+1$.

We put $z_{m(n)+1}=y$ and enumerate $z_{m(n)+1},..., z_{m(n+1)}$  the set $A_{n+1}A_{n+1}^{-1}\setminus (A_{n+1}\bigcup \{z_{0}, ..., $ $ z_{m(n)} \} ) $. Then we choose closed in $H$ neighborhoods $V_{m(n)+1}, ...,   V_{m(n+1)} $ of $e$ such that $z_{i}V_{i}\bigcap $ $A_{n+1} = \emptyset$, $\mu(V_{i})< r_{i}$  for each $i\in\{m(n)+1, ..., m(n)\}$. Thus,  (4), (5)  are satisfied for $k=n+1$.

After  $\omega$ steps we  get the desired $A=\{x_{n}, y_{n}: n\in\omega\}: G= AA^{-1}$ by (1), $A$ is discrete by (2), $A$  is closed by (4).

{\it Case 2.}  We take a neighborhood  $U$  of $e$ such that  $G \neq FU$  for every finite subset $F$  of $G$,  and pick a neighborhood $V$  of $e$  such that $VV^{-1} \subset U$.  Then we choose inductively a sequence $(x_{n})_{n\in\omega}$  in  $G$  such that $\{x_{n}, g_{n}x_{n}\}V\bigcap \{x_{i}, g_{i}x_{i}\}V=\emptyset$ for each  $i< n$.  The set  $A= \{x_{n}, g_{n}x_{n}: n\in\omega\}$  is closed, discrete and  $G=AA^{-1}$.
\vspace{3 mm}

{\bf Question.} {\it Can every countable topological group $G$  be factorized into two (close) discrete subsets $A$  and $B$: $G=AB$  and the subsets $\{aB: a\in A\}$  are pairwise disjoint?}
\vspace{3 mm}

{\bf Remark.}  A subset  $A$  of a group $G$  is called   {\it thin} if $gA \bigcap A$  is finite for each $g\in  G\setminus  \{e\}$. Can every countable topological group be generated by some thin closed discrete subset [2, Question 5]? In both cases of the proof, we have infinitely many possibilities to prolong $A_{n}$ to $A _{n+1}$, so the set $A$  in Theorem can be chosen thin.
\vspace{3 mm}

\centerline{\bf References}
\vspace{3 mm}

[1]  W. Comfort, S. Morris, D. Robbie, S. Svetlichny, M. Tkachenko, {\it Suitable sets for topological groups}, Topology Appl. {\bf 86} (1998)  25-46.

[2]  I. Protasov, {\it Thin subset of topological groups}, Topology Appl.  {\bf 160} (2013),  1083-1087.

\vspace{6 mm}
Department of Cybernetics, Kiev University.

Prospect Glushkova 2, corp. 6,

03680 Kyiv, Ukraine

e-mail:  I.V. Protasov@gmil.com

\end{document}